\documentclass{article}
\usepackage{amsmath,amsthm,amsfonts,amssymb}
\usepackage[all]{xypic}

\theoremstyle{plain}
\newtheorem{prop}{Proposition}
\newtheorem{thm}{Theorem}
\newtheorem{cor}{Corollary}
\newtheorem{lem}{Lemma}

\theoremstyle{definition}
\newtheorem{remark}{Remark}
\newtheorem{example}{Example}
\newtheorem{defi}{Definition}

\newcommand{\vsp}{\vspace{0.05in}}

\begin{document}

\title{Higher Derived Brackets and Deformation Theory I}
\author{F\"{u}sun Akman and Lucian M. Ionescu \\ {\small Illinois State University} }

\maketitle

\begin{abstract}
The existing constructions of derived Lie and sh-Lie brackets involve multilinear maps that are used to define higher order differential operators. In this paper, we prove the equivalence of three different definitions of higher order operators. We then introduce a unifying theme for building derived brackets and show that two prevalent derived Lie bracket constructions are equivalent. Two basic methods of constructing derived strict sh-Lie brackets are also shown to be essentially the same. So far, each of these derived brackets is defined on an abelian subalgebra of a Lie algebra. We describe, as an alternative, a cohomological construction of derived sh-Lie brackets. Namely, we prove that a differential algebra with a graded homotopy commutative and associative product and an odd, square-zero operator (that commutes with the differential) gives rise to an sh-Lie structure on the cohomology via derived brackets. The method is in particular applicable to differential vertex operator algebras.
\end{abstract}

\section{Summary}

The derived bracket constructions of Kosmann-Schwarzbach \cite{KS2,KS} and the higher derived bracket constructions of T.~Voronov \cite{TV,TV2} have several common ingredients. One such ingredient is the use of multilinear maps that define the concept of ``higher order differential operator'' on a commutative associative algebra. There are three main definitions of higher order operators (Definitions ~\ref{gr}, \ref{ko}, and \ref{ak}), and we give a proof of the equivalence of these definitions in the ``Equivalence Theorem'', which is our first main result. We acknowledge that the statement is familiar to experts, yet our present proof is quite possibly the first complete one. 

The second common ingredient in derived bracket constructions is the decomposition of a Lie algebra, such as the Lie algebra of endomorphisms of a unital algebra, into a linear direct sum. We state in Lemma~\ref{decompofendtwo} that the endomorphism algebra $\mbox{End}(A)$ of any unital left pre-Lie algebra $A$ has a decomposition where one summand -identified with $A$- is the Lie subalgebra $\ell_A$ of left multiplication maps and the other is the space of all maps annihilating unity. At this point, the {\it universal theme} of existing constructions emerges: the algebra $A$, usually abelian with respect to the Lie bracket on $\mbox{End}(A)$ or some other similar larger Lie algebra, is equipped with derived brackets that are obtained by modifying the original bracket by a derivation $d$ (namely, $[da,b]$). In particular, the derived Lie bracket constructions summarized in Theorem~\ref{above} involve either a semidirect sum decomposition or the adjoint of a second-order operator as the special derivation (Propositions~\ref{kosscone} and \ref{kossctwo}). In the light of the universal theme, we prove the equivalence of the two approaches in Propositions~\ref{kosschone} and \ref{kosschtwo}.

We also review two derived-bracket examples due to Kosmann-Schwarzbach and present two new ones. In Example~\ref{hoch}, we show that the Hochschild complex of an associative algebra $(A,m)$ (where the differential is the adjoint of $m$) produces as a derived bracket the commutator of $m$ on $A$.

T.~Voronov's ``First and Second Higher Derived Brackets Theorems'' (Theorems \ref{vorone} and \ref{vortwo}) for strict sh-Lie algebras are stated next. We again show that the two statements imply each other by changing the underlying algebra structure (Proposition~\ref{lemvorone} and \ref{lemvortwo}). An important example of sh-Lie construction on a graded commutative associative algebra (first observed in~\cite{BDA}) is given in Proposition~\ref{voronovex}. Moreover, a modification of the Hochschild complex example (Example~\ref{hochtwo}) supplies yet another proof of the well-known fact that the symmetrization of a sh-associative structure is sh-Lie (the original proof is in~\cite{LM}).

Unital left pre-Lie algebras as starting points (instead of Lie algebras) abound until the last section, where we lift this restriction, for several reasons. First of all, the existence of a unity is essential in some bracket definitions and the splitting of the endomorphism algebra. Second, the definitions of several brackets require a product whose commutator is Lie, and only then we can make comparisons between various brackets. Third, every example we study is naturally unital left pre-Lie or can be linearly embedded in such an algebra (e.g.~a Lie algebra is embedded into its universal enveloping algebra). We note that the notions of pre-Lie and left pre-Lie do not necessarily coincide in this paper. In Example~\ref{nabla} we produce a type of pre-Lie algebra that is not always left pre-Lie.
\begin{defi}
Let $(A,m)$ be a graded algebra with bilinear product $m$. We will call the product (and the algebra) {\bf pre-Lie} if the commutator 
$ [a,b]=m(a,b)-(-1)^{|a||b|}m(b,a)$
is a Lie bracket. The product is called {\bf left pre-Lie} if the left multiplication operators $\ell_a=m(a,-)\in\mbox{End}(A)$ satisfy the relation
$ [\ell_a,\ell_b]=\ell_a\ell_b-(-1)^{|a||b|}\ell_b\ell_a=\ell_{ab-(-1)^{|a||b|}ba}=\ell_{[a,b]}.$] Similarly, the product is called {\bf right pre-Lie} if the right multiplication operators $r_b=m(-,b)\in\mbox{End}(A)$ satisfy the relation
$ [r_b,r_c]=r_br_c-(-1)^{|b||c|}r_cr_b=r_{cb-(-1)^{|b||c|}bc}=r_{[c,b]}.$
\end{defi}
%
\begin{remark} Left and right pre-Lie algebras $(A,m)$ are mapped onto the vector space of left/right  multiplication operators in $\mbox{End}(A)$ by a Lie homomorphism and a Lie anti-homomorphism respectively. 
If $(A,m,1)$ is unital left or right pre-Lie, then the map onto multiplication operators is also injective, the inverse given by $\ell_a\mapsto \ell_a(1)=a$ or $r_b\mapsto r_b(1)=b$ respectively. 
\end{remark}
Our second main result, namely the {\it Cohomological Derived Brackets Theorem} (Theorem~\ref{CDBT}), gives an alternative construction of a sh-Lie structure on the cohomology $H(V,Q)$ of a differential algebra $(V,m,Q)$. We assume that $m$ descends to a commutative associative product on $H(V,Q)$, and that there exists an odd, square-zero derivation $\Delta$ of $m$ on $V$ (the differential of the derived bracket). We now eliminate the condition ``unital pre-Lie'', although topological vertex operator algebras (TVOA's) are important examples of unital left pre-Lie algebras to which the new construction can be applied. The cohomological construction partially answers T.~Voronov's inquiry~\cite{TV,TV2} as to how his brackets could be modified for nonabelian Lie subalgebras; instead of modifying the techniques that were useful in the case of abelian subalgebras, we may want to construct derived brackets on the cohomology of a suitable differential graded pre-Lie algebra.

All vector spaces and algebras will be assumed to be over some field of characteristic zero (e.g.~{\bf C}) for simplicity. The symbol $\Box$ will denote the end of a proof. 
\section{The Equivalence Theorem}
\subsection{Definitions and Statement}

A central language that is used in building derived brackets is that of higher order differential operators on commutative associative algebras. 
Let $(A,m,1)$ denote an algebra with underlying vector space $A$, bilinear multiplication $m$, and two-sided unity $1$. Also let $\Delta$ be a linear map in $\mbox{End}(A)$ (the latter is equipped with the Lie bracket that is the commutator of the composition product). For $a\in A$, we will denote the left multiplication operator in $\mbox{End}(A)$ by $\ell_a$. If the vector space $A$ and hence the Lie algebra $(\mbox{End}(A),[-,-]$) are graded, then we will use the Koszul sign convention. We assume even grading for the rest of this section, but reserve the right to revert to the graded case elsewhere, as all results continue to hold with proper sign modifications.

\begin{defi} If $(A,m,1)$ is commutative and associative, then the linear map $\Delta$ on $A$ is called a {\bf differential operator of order at most $k$} if any linear operator $\Gamma_{\Delta}^r(a_1,a_2,\dots,a_r)$ in $\mbox{End}(A)$ (where $a_i\in A$) defined by
\[ \Gamma_{\Delta}^r(a_1,a_2,\dots,a_r)=[\cdots [\,[\Delta,\ell_{a_1}],\ell_{a_2}],\cdots,\ell_{a_r} ]\]
is identically zero for $r\geq k+1$ . (Attributed to Grothendieck~\cite{Gr}) \label{gr}\end{defi}

\begin{remark} A differential operator $T$ of order at most zero is a left multiplication operator (namely, left multiplication by $T(1)$). 
\end{remark}
\begin{remark} Commutativity or associativity is not essential in the definition; $m$ may be taken to be pre-Lie with unity, for example, in certain applications. 
\end{remark}
\begin{remark}
Being of higher order can be defined recursively: any operator of order $\leq -1$ must be zero, and $\Delta$ is of order $\leq k$ if and only if $[\Delta,\ell_a]$ is of order $\leq k-1$ for all $a$. 
\end{remark}
See \cite{Sa} for Sardanashvily's generalization to higher order differential operators in $\mbox{Hom}_{\cal A}(P,Q)$ where ${\cal A}$ is an algebra over a commutative associative ring and $P,Q$ are ${\cal A}$-bimodules. Now define $\lambda:A\rightarrow A\otimes A$ as $\lambda(a)=(1\otimes a-a\otimes 1)$ and extend it multiplicatively to $\lambda^r:A^{\otimes r}\rightarrow A\otimes A$ using the product in $A\otimes A$. 

\begin{defi} If $(A,m,1)$ is commutative and associative, then the linear map $\Delta$ on $A$ is called a {\bf differential operator of order at most $k$} if any vector $\Psi_{\Delta}^r(a_1,a_2,\dots,a_r)$ in $A$ (where $a_i\in A$) defined by
\[ \Psi_{\Delta}^r(a_1,a_2,\dots,a_r)=m\circ (\Delta\otimes \mbox{id}_A)\lambda^r(a_1\otimes\cdots\otimes a_r)
\]
is  zero for $r\geq k+1$ .
(Attributed to Koszul~\cite{Ko}) \label{ko}
\end{defi}

\begin{remark}
Commutativity is not essential in the definition. Associativity can be circumvented by careful inductive definition of the product in $A$ and $A\otimes A$, e.g. from right to left.\end{remark}

\begin{defi} If $(A,m,1)$ is any unital algebra ($m$ suppressed), then the linear map $\Delta$ on $A$ is called a {\bf differential operator of order at most $k$} if any vector $\Phi_{\Delta}^r(a_1,a_2,\dots,a_r)$ in $A$ (where $a_i\in A$) defined by the recursive formula
\begin{eqnarray}
&&\Phi_{\Delta}^1(a)=\Delta(a)-\Delta(1)\, a\nonumber\\
&& \Phi_{\Delta}^2(a,b)=\Phi_1(ab)-\Phi_1(a)b-a\Phi_1(b)\nonumber\\
&&\vdots\nonumber\\
&& \Phi_{\Delta}^{r+1}(a_1,\dots,a_r,a_{r+1})=\Phi_{\Delta}^{r}(a_1,\dots,a_ra_{r+1})-\Phi_{\Delta}^{r}(a_1,\dots,a_r)a_{r+1}\nonumber\\
&&-a_r\Phi_{\Delta}^{r}(a_1,\dots,a_{r+1})\nonumber\end{eqnarray}
is  zero for $r\geq k+1$ .
(Akman~\cite{Ak3}) \label{ak}
\end{defi}

\begin{remark} This definition is designed to work with noncommutative, nonassociative algebras. All three definitions will be used later in contexts where the arguments $a_i$ are in a commutative subalgebra of an associative/pre-Lie algebra (or an abelian subalgebra of a Lie algebra) but the output may be in the larger algebra.
\end{remark}

\begin{thm}[Equivalence Theorem]\label{equi} If $(A,m,1)$ is a commutative associative unital algebra and $\Delta$ is a linear map on~$A$, then the three definitions of higher order differential operator for $\Delta$ are equivalent. That is, the following holds identically for variables $a_i$ in $A$:
\begin{equation}\label{firstline} \Gamma_{\Delta}^r(a_1,\dots,a_r)1=
\Psi_{\Delta}^r(a_1,\dots,a_r)=
\Phi_{\Delta}^r(a_1,\dots,a_r)\end{equation}
(note evaluation at unity for the $\Gamma$ operator). All three expressions are equal to the sum
\begin{equation} \sum_{\sigma\in S_{k,r-k}}
(-1)^k a_{\sigma(1)}\cdots a_{\sigma(k)}\Delta(a_{\sigma(k+1)}\cdots a_{\sigma(r)}).\label{equ}\end{equation}
\end{thm}
Here $S_{k,r-k}$ is the subset of the symmetric group $S_r$ consisting of the $(k,r-k)$-{\bf unshuffles}. These are by definition permutations $\sigma$ in $S_r$ such that  
$\sigma(1)<\cdots <\sigma(k)$ and $\sigma(k+1)<\cdots <\sigma(r)$ 
where $0\leq k\leq r$.
\begin{remark} The second equality in Equation~(\ref{firstline}) is given in \cite{Ak3,BDA} but not proven; also Bering, Damgaard, and Alfaro mention that the first and third expressions in~(\ref{firstline}) are equivalent in~\cite{BDA} and give a proof for the lowest identities~only. 
\end{remark}
\begin{remark} If the condition of commutativity is removed, then $\Psi$ will still have the form in Eq.~(\ref{equ}). For $\Gamma$ operators, only the order of factors on the left-hand side of $\Delta$ in (\ref{equ}) will be completely reversed. The $\Phi$ operators will not look like either. We do not have a general description of the relationships between the operators in the absence of commutativity and associativity; the definitions would involve arbitrary choices and the formulas would be extremely technical. We have found it easier to compare the operators on a case-by-case basis.
\end{remark}
\subsection{Proof}


\begin{lem} Let $(X,m,1)$ be a unital associative algebra. Define the associative product in $X\otimes X$ by 
$ (x\otimes y)(z\otimes w)= xz\otimes yw$ 
for $x,y,z,w\in X$. Then for all $r\geq 1$ and $x_1,\dots,x_r\in X$, we have
\begin{eqnarray} &&(1\otimes x_1-x_1\otimes 1)\cdots (1\otimes x_r-x_r\otimes 1)\\ &&=\sum_{\sigma\in S_{k,r-k}}
(-1)^k x_{\sigma(1)}\cdots x_{\sigma(k)}\otimes x_{\sigma(k+1)}\cdots x_{\sigma(r)}.\nonumber\end{eqnarray}
 \label{prelim}\end{lem}
{\it Proof.} By simple induction. An empty product is equal to unity. $\Box$
\begin{lem} Let $(X,m,1)$ be a commutative subalgebra of a unital associative algebra $(Y,m,1)$, $x_1,\dots,x_r\in X$, $y\in Y$, and $[-,-]$ be the commutator of $m$ on $Y$. Then we have
\begin{eqnarray} && [\dots [\, [ y,x_1],x_2],\dots ,x_r]=m\circ (\mbox{id}\otimes\ell_y) (1\otimes x_1-x_1\otimes 1)\cdots (1\otimes x_r-x_r\otimes 1)\nonumber\\ &=& \Psi_{\ell_y}(x_1,\dots,x_r).\nonumber\end{eqnarray}
\label{prelimtwo}
\end{lem}
{} {\it Proof.} Induction on $r$. First, we have
\begin{eqnarray} &&m\circ(\mbox{id}\otimes \ell_y) (1\otimes x_1-x_1\otimes 1)\\ &&= m(1\otimes yx_1-x_1\otimes y )\nonumber \\ &&=yx_1-x_1y=[y,x_1].\nonumber \end{eqnarray}
By induction, if the statement holds for some $r\geq 1$, then
\begin{eqnarray}&& [\dots [\, [ y,x_1],x_2],\dots ,x_{r+1}]\\
&=&[\dots [\, [ y,x_1],x_2],\dots ,x_r]{x_{r+1}}-{x_{r+1}}[\dots [\, [ y,x_1],x_2],\dots ,x_r]\nonumber \\
&=&m\circ (\mbox{id}\otimes\ell_y) (1\otimes x_1-x_1\otimes 1)\cdots (1\otimes x_r-x_r\otimes 1)(1\otimes {x_{r+1}})\nonumber\\
&&- m\circ (\mbox{id}\otimes\ell_y) (1\otimes x_1-x_1\otimes 1)\cdots (1\otimes x_r-x_r\otimes 1)  ({x_{r+1}} \otimes 1)  \nonumber\\
&=&m\circ (\mbox{id}\otimes\ell_y) (1\otimes x_1-x_1\otimes 1)\cdots (1\otimes x_{r+1}-x_{r+1}\otimes 1) 
\nonumber \end{eqnarray}
(multiplication by $x_{r+1}$ is on the left or right of the tensor product depending on whether we want it before or after $y$). $\Box$ 

\begin{lem} Let the hypotheses be as in the last Lemma, and replace $\ell_y$ with a generic linear operator $\Delta$ on $Y$. Then we have
\begin{eqnarray} \Phi^r_{\Delta}(x_1,\dots ,x_r)&=&m\circ(\mbox{id}\otimes \Delta) (1\otimes x_1-x_1\otimes 1)\cdots (1\otimes x_r-x_r\otimes 1)\nonumber\\
&=& \Psi_{\Delta}^r(x_1,\dots,x_r), \end{eqnarray}
where $\Psi^r_{\Delta}$ is as in Definition~\ref{ko} and $\Phi^r_{\Delta}$ is as in Definition~\ref{ak}.\label{prelimthree}
\end{lem}
%
%
{} {\it Proof.} Let $r=1$. Then 
\[ m\circ (\mbox{id}\otimes \Delta) (1\otimes x_1-x_1\otimes 1)=m(1\otimes \Delta(x_1)-x_1\otimes\Delta(1))=\Delta(x_1)-x_1\Delta(1)=\Phi_{\Delta}^1(x_1).\]
Next, assume that the statement holds for some $r\geq 1$. We have
\begin{eqnarray}&& \Phi_{\Delta}^{r+1}(x_1,\dots,x_r,x_{r+1})\\
&=&\Phi_{\Delta}^r(x_1,\dots,x_rx_{r+1})-\Phi_{\Delta}^r(x_1,\dots,x_r)x_{r+1}- x_r\Phi_{\Delta}^r(x_1,\dots,x_{r+1})
\nonumber \\
&=&m\circ(\mbox{id}\otimes\Delta) (1\otimes x_1-x_1\otimes 1)\cdots (1\otimes x_rx_{r+1}-x_rx_{r+1}\otimes 1)\nonumber\\
&&- m\circ(\mbox{id}\otimes\Delta) (1\otimes x_1-x_1\otimes 1)\cdots (1\otimes x_r-x_r\otimes 1)(x_{r+1}\otimes 1)\nonumber\\
&&- m\circ(\mbox{id}\otimes\Delta) (1\otimes x_1-x_1\otimes 1)\cdots (1\otimes x_{r+1}-x_{r+1}\otimes 1)(x_r\otimes 1)\nonumber\\
&=&m\circ(\mbox{id}\otimes\Delta) (1\otimes x_1-x_1\otimes 1)\cdots (1\otimes x_{r-1}-x_{r-1}\otimes 1)[1\otimes x_rx_{r+1}- x_rx_{r+1}\otimes 1\nonumber\\
&&- (1\otimes x_r-x_r\otimes 1)(x_{r+1}\otimes 1)- (1\otimes x_{r+1}-x_{r+1}\otimes 1)(x_r\otimes 1)],\nonumber\end{eqnarray}
where the expression within the square brackets is precisely
\[ (1\otimes x_{r}-x_{r}\otimes 1)(1\otimes x_{r+1}-x_{r+1}\otimes 1).\;\;\;\Box\] 
{\it Proof of the Equivalence Theorem.} By Lemma~\ref{prelim}, the operators $\Psi^r_{\Delta}$ are of the form given in Eq.~(\ref{equ}); the map $\Delta$ can be placed in front of either grouping. The equivalence of Definitions \ref{ko} and \ref{ak} is given by Lemma~\ref{prelimthree}. Finally, by Lemma~\ref{prelimtwo}, the operators $\Gamma^r_{\Delta}$ satisfy 
\begin{equation} \Gamma^r_{\Delta}(a_1,\dots ,a_r)=M\circ(\mbox{Id}\otimes\ L_{\Delta}) (\mbox{id}\otimes \ell_{a_1}-\ell_{a_1}\otimes\mbox{id})\cdots (\mbox{id}\otimes \ell_{a_r}-\ell_{a_r}\otimes\mbox{id}) \label{gamop}\end{equation}
where the product $m$ is now composition in $\mbox{End}(A)$ (commutative on left multiplication operators), $M$ is the product in $\mbox{End}(A)\otimes\mbox{End}(A)$ described in Lemma~\ref{prelim}, $\mbox{id}=\ell_1$, $\mbox{Id}$ is the identity operator on (not in!) $\mbox{End}(A)$, and $L_{\Delta}$ is left multiplication by $\Delta$ defined on $\mbox{End}(A)$.
Consequently, Lemma~\ref{prelim} is applied to prove our case for the $\Gamma^r_{\Delta}$. $\Box$


\subsection{Corollaries}

\begin{cor} The Equivalence Theorem holds even when the arguments are in a commutative subalgebra $A$ of a larger unital associative algebra $B$ and $\Delta\in\mbox{\em End}(B)$.\label{corol}
\end{cor}
%
%
Let $\mbox{Diff}^{\, 0}(A)=\ell_A$, $\mbox{Diff}^{\, 1}(A)$, $\mbox{Diff}^{\, 2}(A)$, ... denote the subspaces of differential operators in $\mbox{End}(A)$ for $(A,m,1)$ as above that consist of operators of orders less than or equal to zero, one, two, ... respectively. 
\begin{cor} For commutative associative $(A,m,1)$, we have the filtration
\[ \mbox{\em Diff}^{\, 0}A\subset \mbox{{\em Diff}}^{\, 1}A\subset\cdots\subset\mbox{\em Diff}^{\, r}A\subset\mbox{\em Diff}^{\, r+1}A\subset\cdots \subset \bigcup_r \mbox{\em Diff}^{\, r}A \subset\mbox{\em End}A \]
under any definition of ``order'' of a differential operator.
\end{cor}
{\it Proof.} By Definitions~\ref{gr} and \ref{ak}, $\Gamma$ and $\Phi$ operators are defined recursively. $\Box$
\begin{cor} Compositions of higher order differential operators on a commutative associative algebra $A$ (under any definition) preserve order:\label{compos}
\[ \mbox{\em Diff}^{\, k}A\circ \mbox{\em Diff}^{\, l}A\subset\mbox{\em Diff}^{\, k+l}A.\]
On the other hand, commutators reduce the total order by one:
\[ [\mbox{\em Diff}^{\, k}A,\mbox{\em Diff}^{\, l}A]\subset\mbox{\em Diff}^{\, k+l-1}A.\]
\end{cor}
{\it Proof.} For compositions we use the following identity for $f,g,h$ in $\mbox{End}(A)$: 
\[ [f\circ g,h]=[f,h]\circ g + f\circ [g,h].\]
The proof is by induction on the total degree and uses Definition~\ref{gr}. 
The proof for brackets is similar and is based on the identity
\[ [\,[f,g],h]=[\,[f,h],g]+[f,[g,h]\,].\;\;\;\Box\]
\begin{remark} \label{remder}
Any element of $\mbox{Diff}^{\, 1}(A)$ for graded commutative associative $A$ with unity can be uniquely written as the sum of a left multiplication operator and a {\bf derivation} $d$ that satisfies the product rule  
$ d(ab)=(da)b+(-1)^{|d||a|}a(db)$ 
and annihilates unity (see Lemma~\ref{decompofend} below). We will make a distinction between the Lie algebra $\mbox{Der}(A)$ of derivations and the algebra $\mbox{Diff}^{\, 1}(A)$ of differential operators of order at most one. 
\end{remark}
Here is another higher-bracket construction related to the $\Gamma$, $\Psi$, and $\Phi$ operators. We define, following Kosmann-Schwarzbach and Voronov, Grothendieck-like operators $B^r_{\Delta}$ as follows:
\begin{defi} We define higher derived brackets
\[ B^r_{\Delta}(a_1,\dots,a_r)=[\dots [\,[\Delta,a_1],a_2],\dots,a_r]\]
where $(L,[-,-])$ is a Lie algebra, $\Delta,a_1,\dots,a_r\in L$, and the $a_i$'s belong to a (possibly abelian) subalgebra $L_0$ of $L$. The outcome need not fall in $L_0$.\label{def4}
\end{defi}
%
%
\begin{lem}\label{lemfour} For $y,x_1,\dots,x_r\in L$, a Lie algebra, such that $x_1,\dots,x_r$ are in an abelian subalgebra, we have
\begin{eqnarray} B^r_{y}(x_1,\dots,x_r)
&=& [\dots [\,[y,x_1],x_2],\dots,x_r]\nonumber\\
&=&m\circ (\mbox{id}\otimes\ell_y) (1\otimes x_1-x_1\otimes 1)\cdots (1\otimes x_r-x_r\otimes 1)\nonumber\\
&=& \Psi_{\ell_y}^r(x_1,\dots,x_r)\nonumber\end{eqnarray}
in the universal enveloping algebra ${\cal U}L$ of $L$. Here $m$ 
denotes the associative product on ${\cal U}L$ whose commutator is the Lie bracket on $L\subset{\cal U}L$. If $L$ is unital associative, then an additional similar result holds, with $m$ the product inducing the Lie bracket on $L$, and $\ell_y$ denoting left multiplication in $L$ (Lemma~\ref{prelimtwo}).
\end{lem}
%
%
\begin{defi} Let $(L,[-,-])$ be a Lie algebra, $a_1,\dots,a_r\in L$ belong to a (possibly abelian) subalgebra of $L$, and $d$ be a linear operator on $L$. Then we define derived brackets
\[ C^r_{d}(a_1,\dots,a_r)=[\dots [da_1,a_2],\dots,a_r].\]
Once again, the image of $C^r_{d}$ need not be in the (abelian) subalgebra.\label{cousin}
\end{defi}
%
%
\begin{lem} Let the hypotheses be as in the last Lemma, and replace $\mbox{\em ad}(y)$ with a linear map $T$ on ${\cal U}L$ that restricts to an endomorphism of $L$. Then we have
\begin{eqnarray} C^r_{T}(x_1,\dots ,x_r)&=&m\circ(\mbox{id}\otimes T) (1\otimes x_1-x_1\otimes 1)\cdots (1\otimes x_r-x_r\otimes 1)\nonumber\\
&=& \Psi_{T}^r(x_1,\dots,x_r),\nonumber
\end{eqnarray}
where $m$ is again the associative multiplication on ${\cal U}L$.
If $L$ is unital associative  as well, then an additional similar result holds, with $m$ the product in $L$ whose commutator is the Lie bracket (Lemma~\ref{prelimthree}).
\label{lemfive}
\end{lem}
\begin{cor}\label{nonewbrackets}
Let $(Y,m,1)$ be unital associative, $(X,m,1)$ a commutative subalgebra, and T be an operator on $Y$. Then for $x_i\in X$ we have 
\begin{eqnarray} && \Phi^r_{T}(x_1,\dots,x_r)\nonumber\\
&=& \Psi^r_T(x_1,\cdots,x_r)\nonumber\\
&=&\Gamma^r_T(x_1,\dots,x_r)1\nonumber\\
&=& B^r_T(\ell_{x_1},\cdots,\ell_{x_r})1\nonumber\\
&=& C^r_{\mbox{\em ad}(T)}(\ell_{x_1},\cdots,\ell_{x_r})1.\nonumber
\end{eqnarray}
\end{cor}



\section{Linear sums of algebras: The universal theme}
\subsection{Semidirect and direct Lie sums}

\begin{defi}If $L_0$ is an ideal and $L_1$ a subalgebra of the Lie algebra $L=L_0\oplus L_1$, consequently satisfying the properties 
$ [L,L_0]\subset L_0$ and $[L_1,L_1]\subset L_1$, 
then we say that $L$ is  a Lie {\bf semidirect sum} (more commonly known as {\bf semidirect product}) of $L_0$ and $L_1$, and denote it by $L=L_0>\!\!\!\lhd\,\, L_1$. \end{defi}
Since $L_1$ acts by derivations on $L_0$, the prototype of a Lie semidirect sum is $A>\!\!\!\lhd\,\,  \mbox{Der}(A)$ where $A$ is an associative, left pre-Lie, or Lie algebra: we define
\[ [a_1+D_1,a_2+D_2]=[a_1,a_2]+D_1(a_2)-D_2(a_1)+[D_1,D_2]\]
for $a_1,a_2\in A$ and derivations $D_1,D_2$. A Lie algebra {\bf direct sum} is a special semidirect sum where both  subalgebras are ideals and the cross-brackets vanish.

A recent analog is the {\bf OCHA's (open-closed homotopy algebras)} in Kajiura and Stasheff \cite{KaS} where one subalgebra is sh-Lie and acts on the other ($A_{\infty}$) by ``derivations''.  This reminds us of the $A_{\infty}$ and $L_{\infty}$ operators that coexist inside a ``weakly homotopy Gerstenhaber algebra'' as defined in \cite{Ak2}.
\vsp


\subsection{Endomorphisms of unital pre-Lie algebras}

\begin{defi} \cite{TV} \label{orderelt}
Let $L_0$ be a subalgebra of a Lie algebra $L$. Then we say that the {\bf order} of an element $\Delta$ of $L$ with respect to $L_0$ is at most $r$ if all expressions
\[ B_{\Delta}^{r+1}(a_1,\dots,a_{r+1})=[\dots [\, [\Delta ,a_1],a_2],\dots ,a_{r+1}]\]
vanish for elements $a_1,\dots,a_{r+1}$ of $L_0$. 
\end{defi}
For a graded commutative associative algebra $A$ embedded into $\mbox{End}(A)$ as an abelian subalgebra (identified with left multiplication operators), the elements of $\mbox{End}(A)$ of order at most $r$ with respect to $\ell_A$ are exactly the differential operators on $A$ of order at most $r$. The subspace of linear endomorphisms of any associative algebra $A$ that commute with all the left and right multiplication operators is called the {\bf centroid} of $A$. 
\begin{lem} Let $(A,m,1)$ be a unital graded commutative associative algebra. Then the graded Lie algebra $\mbox{\em End}(A)$ is the linear direct sum of the abelian subalgebra $\ell_A$ that consists of left multiplications by elements of $A$ and the subalgebra $\mbox{\em Ann}(1)$ of linear maps on $A$ that annihilate 1:
\[ \mbox{\em End}(A)=\ell_A\oplus \mbox{\em Ann}(1).\]
The subalgebra $\ell_A$ is the centroid of $A$. Neither $\ell_A$ nor $\mbox{\em Ann}(1)$ is an ideal. \label{decompofend}
\end{lem}
{\it Proof.} Clearly, the projection $P:\mbox{End}(A)\rightarrow A$ given by $T\mapsto T(1)$ induces an isomorphism of $\ell_A$ onto $A$ and has kernel $\mbox{Ann}(1)$, easily seen to be a subalgebra. The rest of the proof is left to the reader.
$\Box$

In fact, $\mbox{End}(A)$ has a similar decomposition for any unital algebra $A$. 
\begin{lem} \label{decompofendtwo} If $(A,m,1)$ be a unital graded algebra, then we have $ \mbox{\em End}(A)=\ell_A\oplus \mbox{\em Ann}(1).$ 
The subspace $\ell_A$ is at the same time the space of differential operators on $A$ of order at most zero under Definition~\ref{ak}. Any differential operator of order at most one can be uniquely written as the sum of a left multiplication operator and a derivation. If $m$ is left pre-Lie or associative, then $\ell_A$ is also a Lie subalgebra of $\mbox{\em End}(A)$ identified with $A$; in this case, we have  
$ \mbox{\em Diff}^{\, 1}(A)=\ell_A>\!\!\!\lhd\,\,  \mbox{\em Der}(A)$.
\end{lem}
\begin{cor}\label{decompofendthree}
Let $L$ be a Lie algebra. Then $L$ is embedded into $\mbox{\em End}({\cal U}(L))$ as the Lie subalgebra of left multiplication operators $\ell_{L}$. If $L_0\subset L$ is a Lie subalgebra, then $\ell_{L_0}$ is a Lie subalgebra of $\ell_{L}$ as well as of $\mbox{\em End}({\cal U}(L))$. 
\end{cor}


\subsection{The Universal Theme}

The {\it universal theme} underlying many derived bracket constructions is as follows: we embed the space $A$ (or~$L_0$) that we are trying to endow with a bracket as an abelian Lie subalgebra into $\mbox{End}(A)$ or some similar Lie algebra. Then we use a (possibly inner) derivation $d$ in the larger space to define the derived bracket on $A$ by $C_d^2(a,b)=[da,b]$ or $B_{\Delta}^2(a,b)=[\,[\Delta,a],b]$. The derived Lie constructions usually make use of a semidirect sum, or bracketing with a second order differential operator, resulting in a bona fide Lie algebra structure on $A$, whereas derived sh-Lie constructions restrict a series of brackets $C^r_d(a_1,\dots,a_r)$ or $B^r_{\Delta}(a_1,\dots,a_r)$ back to $A$ by using a projection onto $A$.


\section{Derived Lie and sh-Lie brackets}


\subsection{Sh-Lie algebras}

All vector spaces and maps are super-graded (see~\cite{TV} for the sign convention). 

\begin{defi} \cite{LS} A (strict) {\bf strongly homotopy Lie algebra} (or {\bf sh-Lie algebra}, a.k.a. {\bf $L_{\infty}$ algebra}) is a super graded vector space $A$ and a sequence of odd graded-symmetric brackets 
$ [-]_1,[-,-]_2,\dots,[-,\cdots,-]_n,\dots $ 
satisfying the relations
\[ \sum_{\sigma\in S_{k,r-k}}(-1)^{\epsilon}[\, [a_{\sigma(1)},\dots ,a_{\sigma(k)}]_k,a_{\sigma(k+1)},\dots ,a_{\sigma(r)}]_{r-k+1}=0                  \]
for each $r\geq 1$; here the {\bf Koszul sign} $(-1)^{\epsilon}$ is given by the product of all factors 
$ (-1)^{|a_{\sigma(i)}||a_{\sigma(j)}|}$ 
for which $i<j$ but $\sigma(i)$ is to the right of $\sigma(j)$. We call the above identity the {\bf Jacobi identity} for sh-Lie algebras.
In particular, the linear map $[-]_1$ is a square-zero, odd derivation of the bilinear bracket $[-,-]_2$. 
\label{linf}
\end{defi}


\subsection{Derived Lie brackets.} \label{derkos}
\subsubsection{Semidirect sums}

The following results have appeared in \cite{KS2,KS}. A (left) {\bf Leibniz algebra} is a vector space with a bilinear bracket for which the (left) adjoint of any element acts as a derivation of the bracket. A Leibniz algebra with an anti-symmetric bracket is a Lie algebra.
\begin{defi} If $(L,[-,-])$ is a graded Lie or Leibniz algebra with a bracket of degree $n$ and an odd, square-zero derivation $d$, then the {\bf derived bracket} on $L$ induced by $d$ is defined by
\[ [a,b]_{d}=(-1)^{n+|a|+1}[da,b]=(-1)^{n+|a|+1}C_d^2(a,b),\]
where $|a|$ is the degree of $a$ (Kosmann-Schwarzbach).\label{kosschdef}
\end{defi}
The parity of the original bracket is reversed in this construction. 
In order to obtain a genuine Lie bracket on a subalgebra $L_0$ of a Lie algebra $L$ (assumed to be abelian under the original bracket $[-,-]$), we need to put the restriction $[dL_0,L_0]\subset L_0$. 
\begin{thm}[Derived Lie Brackets Theorem] \cite{KS2,KS} \label{above}

(i) If $(L,[-,-])$ is a Lie or Leibniz algebra as above, then the derived bracket on $L$ induced by $d$ satisfies the (left) Leibniz property.

(ii) Such a derivation $d$ of $(L,[-,-])$ is also a derivation of the derived bracket.

(iii) Let $L_0$ be an abelian subalgebra of a Lie algebra $(L,[-,-])$, and $d$ be an odd, square zero derivation of $L$ such that $[dL_0,L_0]\subset L_0$. Then the restriction of the derived bracket to $L_0$ is  graded symmetric, and we obtain a graded Lie algebra.

(iv) Another way to obtain a Lie bracket is to pass to the quotient of $L$ by $dL$.

\end{thm}
The notion of a derived bracket arose in Kosmann-Schwarzbach's work in the following form (see \cite{KS2,KS}): given a graded vector space $L$ and a linear map $f:L\rightarrow\mbox{End}(L)$, we have the multiplication 
$ [a,b]_f=f(a)b$ 
induced by $f$ on $L$. 
(Although $f(a)$ corresponds to the left multiplication operator $\ell_a$, it is sometimes taken as ad($a$).) We will instead consider a Lie {\it embedding} of an abelian Lie algebra $L_0$ into a larger Lie algebra $L$, possibly $\mbox{End}(L_0)$. We will assume that a particular linear complement $L_1$ of $L_0$ in $L$ is given, but the existence of $L_1$ is not technically necessary for Kosmann-Schwarzbach's constructions of the type we discuss in this subsection. In most examples, such $L_1$ exist as subalgebras. Then we claim that the essence of the specific constructions in \cite{KS} is a linear direct sum $L=L_0\oplus L_1$ where $L$ is a Lie algebra and $L_0$ is an abelian subalgebra. In particular, we have
\begin{prop}\label{kosscone} Let $L_0$ and $L_1$ be Lie algebras ($L_0$ abelian).
If $L=L_0>\!\!\!\lhd\,\, L_1$ is a Lie semidirect sum with an odd, square-zero derivation $d$ of $L$, then the derived bracket $[-,-]_d$ on $L$ naturally restricts to $L_0$ and becomes a Lie bracket. In case $L_1$ is also an ideal and the sum is direct, the derived bracket vanishes on $L_0$.
\end{prop}
\begin{example} \cite{KS} Let $\Omega^{\bullet}(M)$ denote the de Rham complex of differential forms on a manifold $M$. 
Cartan's formulas
\[ [d,d]=0, \;\;\; [\iota_X,\iota_Y]=0,\;\;\; {\cal L}_X=[\iota_X,d],\;\;\; [{\cal L}_X,d]=0,\;\;\; [{\cal L}_X,\iota_Y]=\iota_{[X,Y]},\]
with the de Rham differential $d$, vector fields $X$, $Y$, substitution operators $\iota_X$, $\iota_Y$, and the Lie derivation ${\cal L}_X$, show us that the Lie bracket $[-,-]$ of vector fields is a derived bracket induced by $d$: 
\begin{equation} \iota_{[X,Y]}=[\,[d,\iota_X],\iota_Y].\label{derham}\end{equation} 
The universal theme appears as follows: each element of the graded Lie subalgebra $\mbox{Der}(\Omega^{\bullet}(M))\subset\mbox{End}(\Omega^{\bullet}(M))$ of degree $k$ is uniquely a sum of the form
\[ D={\cal L}_K+\iota_L,\quad K\in \Omega^{k}(M;TM),\quad L\in \Omega^{k+1}(M;TM)                \]
(e.g.~\cite{KMS}). That is,
\[ \mbox{Der}(\Omega^{\bullet}(M))=\Omega^{\bullet}(M;TM)\oplus \Omega^{\bullet +1}(M;TM)={\cal L}_M\oplus \iota_M,\]
where the first summand is a subalgebra and the second is an ideal. In Equation~(\ref{derham}) we have $d$ as an odd, square zero derivation and the $\iota_X$'s forming an abelian subalgebra of $\iota_M$, where vector fields $X\in\mbox{Vect}(M)$ are identified with $\iota_X$. Therefore, under this identification, the bracket $[-,-]$ on the Lie algebra $\mbox{Vect}(M)$ does come from a derived bracket on a larger Lie algebra -a semidirect sum- where $\mbox{Vect}(M)$ is an odd abelian subalgebra. \label{ksex}
\end{example}

\begin{example} \label{nabla}
Let the differential geometric notation be as before, including $L=\mbox{Vect}(M)$. Every linear connection $\nabla$ is a multiplication-generating map $f$ in the Kosmann-Schwarzbach sense. That is, 
$ \nabla :L\rightarrow \mbox{End}(L)$, 
where the induced multiplication $XY=\nabla_X(Y)$ on $L$ is pre-Lie, but in general not left or right pre-Lie, for a torsion-free or symmetric connection (we hereby fix such a connection $\nabla$): the condition
\[ T(X,Y)=\nabla_X(Y)-\nabla_Y(X)-[X,Y]=0\]
shows us that the commutator of $XY$ on $L$ is the Lie bracket of vector fields on $L$. We then examine the curvature tensor
$K(X,Y)=[\nabla_X,\nabla_Y]-\nabla_{[X,Y]},\label{curv}$  
which defines a bilinear map on $L$ with values in $\mbox{End}(L)$. If $\nabla$ also has zero curvature, then the Lie bracket in the subspace $\nabla_L$ of $\mbox{End}(L)$ closes in $\nabla_L$ and defines a homomorphism of $L$ as a Lie algebra into $\mbox{End}(L)$. The Jacobi identity in $\nabla_L$ is also known as the {\bf Bianchi identity}. 

Note that in the case of zero torsion and zero curvature, the product $XY=\nabla_X(Y)$ on $L$ is in fact precisely left pre-Lie, because 
\begin{eqnarray} &&\{ \,(XY)Z-X(YZ)\,\} -\{ \, (YX)Z-Y(XZ)\,\}\nonumber\\
&=&(\nabla_{\nabla_X(Y)}(Z)-\nabla_X\nabla_Y(Z))-(\nabla_{\nabla_Y(X)}(Z)-\nabla_Y\nabla_X(Z))\nonumber\\
&=& \nabla_{\nabla_X(Y)-\nabla_Y(X)}(Z)-[\nabla_X,\nabla_Y](Z)\nonumber\\
&=& \nabla_{[X,Y]}(Z)-\nabla_{[X,Y]}(Z)=0.\nonumber\end{eqnarray}

Applying the last example above to our case, we deduce that the bracket on the homomorphic image of $L=\mbox{Vect}(M)$ is inherited from a derived bracket on $\mbox{Der}(\Lambda^{\bullet}(\mbox{End}(L)'))$, thanks to the connection $\nabla$. We cannot say that $L\cong\nabla_L$, though, since there is no unity with respect to the product $\nabla_XY$ (also see~\cite{Ion1,Ion2}). 

\end{example}

\subsubsection{BV constructions}

We recall that a {\bf Gerstenhaber (G) algebra} is a vector space with a graded commutative associative algebra structure and an additional odd bracket (the G-bracket) that makes it into a graded Lie algebra. The following condition also holds: bracketing with an element of this space is a derivation of the commutative associative product (then a G-algebra is a graded version of a Poisson algebra). We also recall that a {\bf Batalin-Vilkovisky (BV) algebra} is a G-algebra where the bracket measures the deviation of an odd, square-zero, second-order operator from being a derivation of the commutative associative product. 

Once again, we are assuming that a particular complement $L_1$ of $L_0$ in the Lie algebra $L$ is given for comparison purposes, but the existence of such a subspace is not technically necessary. 
Here is another way of ensuring $[dL_0,L_0]\subset L_0$: Let $L_0$ be an abelian subalgebra of $L$ such that {\it the centralizer $Z_{L}(L_0)$ of $L_0$ is itself}. 
Choose the bracket-generating derivation $d$ to be the adjoint of an odd, square-zero element $\Delta$ of $L$ that is of order at most two with respect to $L_0$. Then by Definition~\ref{orderelt}, we obtain the characterization
\[ \mbox{$\Delta$ is of order}\leq 2
\Leftrightarrow [\, [\, [\Delta,L_0],L_0],L_0]=0
\Leftrightarrow [\, [\Delta,L_0],L_0]\subset Z_{L}(L_0)=L_0.\]
\begin{prop} \label{kossctwo}
If $L_0$ is an abelian subalgebra of a Lie algebra $L$ which is its own centralizer, and $\Delta\in L$ is odd, square-zero element, and of order at most two with respect to $L_0$, then the derived bracket $[-,-]_{\mbox{\em ad}(\Delta)}$ on $L$ restricts to $L_0$ and becomes a Lie bracket. In case $L$ is a Lie direct sum of ideals, the derived bracket vanishes on $L_0$.
\end{prop}
\begin{remark} \label{fakecon}
Even if the condition $Z_L(L_0)=L_0$ is not satisfied, we will label a derived bracket as ``BV-type'' as long as we have a $\Delta$ of order at most two with $[\, [\Delta,L_0],L_0]\subset L_0$. 
\end{remark}
We have already seen an example of self-centralizing $L_0$, namely a unital graded commutative associative algebra $A$, identified with $\ell_A$ inside $\mbox{End}(A)$:
\begin{example} \cite{KS} 
Let $(A,m,1)$ be a graded commutative associative algebra. An inner derivation via an odd, square-zero, linear map $\Delta\in\mbox{End}(A)$ determines a derived bracket $[-,-]_{\Delta}$ on $\mbox{End}(A)$ that is Leibniz; if $\Delta$ is also a second-order differential operator, then the restriction of $[-,-]_{\Delta}$ to $\ell_A=A$ becomes {\it the} BV bracket, because it measures the deviation of $\Delta$ from being a derivation of $m$ (also see~\cite{Ak3}): we have
\[ [\ell_a,\ell_b]_{\Delta}1=[\,[\Delta,\ell_a],\ell_b]1=
\Gamma_{\Delta}^2(a,b)1=\Phi_{\Delta}^2(a,b) \]
by the Equivalence Theorem, and the last expression is the definition of the BV bracket in~$A$. Applying the last Proposition to $\mbox{End}(A)=\ell_A\oplus\mbox{Ann}(1)$, we obtain a derived bracket on $A$.
\label{ksex}
\end{example}

Here is a new example where restriction to $L_0$ is due to a favorable~grading:

\begin{example} The Lie bracket on an associative algebra $(A,m)$ is a derived bracket induced from its Hochschild complex $ C^{\bullet}(A)=\bigoplus_{n\geq 0}\mbox{Hom}(A^{\otimes n},A)$ with the pre-Lie composition and the Gerstenhaber bracket. The subspace $\mbox{Hom}(A^{\otimes 0},A)=A$ is an abelian subalgebra with respect to the G-bracket:
\[ [a,b]=a\circ b-b\circ a=0-0=0\]
(compositions by an element of $A$ on the left vanish by definition). This is a case where higher-order elements with respect to the abelian subalgebra are readily identified.  Maps $f\in C^{1}(A)$ are of order at most one with respect to $A$:
\[ [f,a]=f\circ a -a\circ f=f(a)-0=f(a)\in A \;\;\mbox{and}\;\;  [\, [f,a],b]\in [A,A]=0.\]
Similarly, bilinear maps $x$ are elements of order at most two with respect to 0-linear maps, because it would take two brackets with members of $A$ to send such maps into $A$: $[\,[x,a],b]=x(a,b)- x(b,a)\in A$, and $B^3_{x}\equiv 0$. In general, every $n$-linear map (as well as every $k$-linear map with $k<n$) is an element of order at most $n$. Thus the associative multiplication $m$ is a square-zero element of order at most two. The {\it degree} of any $n$-linear map is by definition $n-1$, so that $m$ is also odd. Now, the natural restriction of the derived bracket $[-,-]_m=[\,[m,-],-]$ to $A$ is the Lie bracket on $A$ associated to $m$: we have
\[ [a,b]_m=[\,[m,a],b]=m(a,b)-m(b,a).\]
Here $A$ is not an ideal as the G-bracket preserves degrees, nor is it true that the centralizer of $A$ is itself; still, we have the desired condition $[\, [m,A],A]\subset A$. In fact, the filtering induced by ``degree'' (which differs from the ``order'' filtering by one) forces the resulting elements to fall into $A$: we have $1+(-1)+(-1)=(-1)$.
 \label{hoch}
\end{example}
We note that the ``BV construction'' follows from the ``semidirect sum'' construction, and vice versa: 
\begin{prop} Let $L_0$ be an abelian subalgebra of some Lie algebra $L$ and $\Delta\in L$ be an odd, square-zero element that is of order at most two with respect to $L_0$. We assume that the derived bracket  
$ [u_0,v_0]_{\Delta}=[\,[\Delta,u_0],v_0]$ 
on $L_0$ obtained by the BV construction is closed on $L_0$ (which may be, for example, due to $Z_L(L_0)=L_0$ or propitious filtering). Then the bracket above is at the same time the derived bracket 
$ [u_0,v_0]_d=[du_0,v_0]$ 
given by the derivation $d=\mbox{\em ad}(\Delta)$ on the semidirect sum 
$ {L^{\ast}}=L_0>\!\!\!\lhd\,\,\mbox{\em ad}[\Delta,L_0]$.
\label{kosschone}
\end{prop}
{\it Proof.} The subspace $[\Delta,L_0]$ of $L$ is easily seen to be a Lie algebra in its own right since $\Delta$ is square-zero. The derivation $d=\mbox{ad}(\Delta)$ sends $L_0$ to $[\Delta,L_0]$ and annihilates the latter, so ${L^{\ast}}$ is closed under it. $\Box$
\begin{prop}\label{kosschtwo}
Let $d$ be an odd, square-zero derivation of the semidirect sum $ L= L_0>\!\!\!\lhd\,\, L_1$ where the ideal $L_0$ is abelian. Then the derived bracket 
$ [u,v]_d=[du,v]$ 
on $L_0$ is also the BV derived bracket
$ [u,v]_{d}=[\,[d,\ell_{u}],\ell_{v}]$ 
where the right-hand side is an expression in $\mbox{\em End}({\cal U}(L))$ that restricts to $L_0$.
\end{prop}
{\it Proof.} The universal enveloping algebra ${\cal U}(L)$
contains $L$ and the abelian subalgebra 
$L_0$ by Corollary~\ref{decompofendthree}.
Any Lie derivation $d$ of $L$ extends uniquely to an associative and Lie derivation of ${\cal U}(L)$.
Then $d\in\mbox{End}({\cal U}(L))$ is still odd, square-zero ($d^2$ and $0$ are two derivations that agree on $L$), and an element of order at most two with respect to $\ell_{L_0}\subset\mbox{End}({\cal U}(L))$, since we have
\[ [\,[d,\ell_{u}],\ell_{v}]=[\ell_{du},\ell_{v}]= \ell_{[du,v]}\in\ell_{L_0}\]
for $u,v\in L_0$. 
$\Box$
%
%


\subsection{Derived sh-Lie brackets}

T.~Voronov's higher derived brackets are defined in \cite{TV} and \cite{TV2} in terms of a projection and form an $L_{\infty}$ algebra rather than a Lie or Leibniz algebra. 
\begin{thm}[First Higher Derived Brackets Theorem: Abelian Case, Inner Derivation]\label{vorone}
Let $L=L_0\oplus L_1$ be a Lie algebra that is a direct sum of two subalgebras as a vector space, and $P$ be the canonical projection onto $L_0$ with kernel $L_1$. Assume moreover that $L_0$ is abelian and $\Delta\in L_1$ is odd and square-zero. Then the operators $P\circ B^r_{\Delta}$ ($r\geq 1$) form an $L_{\infty}$ algebra structure on $L_0$. 
\end{thm}
\begin{remark}
The condition $\Delta\in L_1$ ensures that the 0-ary bracket, $P(\Delta)$, vanishes. In other words, we have a ``strict'' $L_{\infty}$ algebra. The original theorem in \cite{TV} involves a general $\Delta\in L$ which may not satisfy this condition. 
\end{remark}
\begin{thm}[Second Higher Derived Brackets Theorem: Abelian Case, Arbitrary Derivation]\label{vortwo}
Let $L=L_0\oplus L_1$ be a Lie algebra that is a direct sum of two subalgebras as a vector space, and $P$ be the canonical projection onto $L_0$ with kernel $L_1$. Assume moreover that $L_0$ is abelian and there exists a linear operator $d$ which is an odd, square-zero derivation of $L$, not necessarily inner, preserving the subalgebra $L_1$.
Then the operators $P\circ C^r_{d}$ ($r\geq 1$) form an $L_{\infty}$ algebra structure on $L_0$.\end{thm}
%
%
\begin{prop}
Theorem~\ref{vorone} implies Theorem~\ref{vortwo}. \label{lemvorone}
\end{prop}
{} {\it Proof.} Assume Theorem~\ref{vorone}. Let $L=L_0\oplus L_1$ be a Lie algebra that is a linear direct sum of two subalgebras, and $P$ be the canonical projection onto $L_0$ with kernel $L_1$. Assume moreover that $L_0$ is abelian and $d\in\mbox{Der}(L)$ is odd and square-zero, preserving the subalgebra $L_1$. We form the Lie algebra
\[ {L^{\ast}}=L>\!\!\!\lhd\,\,\mbox{\bf C}d=L_0\oplus L_1\oplus \mbox{\bf C}d\]
where $[x,d]=-d(x)$ for $x\in L$ by definition of semidirect sum and $[d,d]=0$. Let ${L_1^{\ast}}=L_1\oplus \mbox{\bf C}d$. Then we have 
$ {L^{\ast}}=L_0\oplus{L_1^{\ast}}$ 
where ${L_1^{\ast}}$ is also a subalgebra (d preserves $L_1$); denote the canonical projection from ${L^{\ast}}$ onto ${L_0}$ by ${P^{\ast}}$. Now $d\in\mbox{Der}(L)$ whereas $\Delta =d\in {L_1^{\ast}}$. Theorem~\ref{vorone} is then applicable to ${L^{\ast}}$, giving us the following $L_{\infty}$ structure on $L_0$:
\begin{eqnarray} && {P^{\ast}}B_{\Delta}^r(a_1,\dots,a_r)\nonumber\\
&=& {P^{\ast}}[\cdots [\, [\Delta,a_1],a_2],\cdots ]\nonumber\\
&=& P [\cdots [da_1,a_2],\cdots ]\nonumber\\
&=& PC_d^r(a_1,\dots,a_r).\;\;\;\Box\nonumber\end{eqnarray}
%
%
%
\begin{prop} \label{lemvortwo}
Theorem~\ref{vortwo} implies Theorem~\ref{vorone}.
\end{prop}
\begin{example}\label{hochtwo}
In Example~\ref{hoch} we saw that any $n$-linear map in the Hochschild-Gerstenhaber Lie algebra of an associative algebra $(A,m)$ would constitute an element of order at most $n$ with respect to $A$. Then the decomposition 
$ C^{\bullet}(A)=A\oplus \bigoplus_{n\geq 1}\mbox{Hom}(A^{\otimes n},A)$ 
helps us define a truncated  $L_{\infty}$ algebra on $A$ for any odd, square-zero multilinear map $\Delta\in C^{\bullet}(A)$. If we allow infinite sums of the form 
$ m=m_1+m_2+\cdots +m_n+\cdots$ 
where $m_n$ is an $n$-linear map, then any given $A_{\infty}$ structure $m$ on $A$ with $[m,m]=0$ (under the {\it suspended degree} convention, with all $m_n$ odd) would induce an $L_{\infty}$ structure on $A$, defined by 
\[ PB_m^r(a_1,\dots,a_r)=P[\cdots [\, [m,a_1],a_2],\cdots ,a_r]=[\cdots [\, [m_r,a_1],a_2],\cdots ,a_r].\]
Clearly, this is the same old $L_{\infty}$ algebra obtained by symmetrizing $m$ (see~\cite{LM}). 
\end{example}
Here is an important $L_{\infty}$ construction of earlier origin:
\begin{prop}[Bering et al.~\cite{BDA}, Kravchenko~\cite{Kr}, T.~Voronov~\cite{TV}] \label{voronovex}
Let $(A,m,1)$ be a graded commutative associative algebra with an odd, square-zero endomorphism $\Delta$. Then the brackets
$\Phi^r_{\Delta}$ form an $L_{\infty}$ structure on $A$.
\end{prop}
\begin{remark} The brackets $\Gamma^r_{\Delta}$ (evaluated at unity) and $\Psi^r_{\Delta}$ would work equally well by the Equivalence Theorem. 
\end{remark}
\begin{remark}
The similar construction in Example~\ref{ksex} ends at the binary bracket because $\Delta$ is of second order. With an  operator of higher order $k$, we can construct $L_k$ algebras ($L_{\infty}$ algebras in which $n$-ary brackets vanish after $n=k$).
\end{remark}
{} {\it Proof of Proposition~\ref{voronovex}.}
We follow Voronov's proof in~\cite{TV}. Let $L=\mbox{End}(A)=\ell_A\oplus\mbox{Ann}(1)$ and define the projection
\[ P:L\rightarrow L,\;\;\; P(T)=\ell_{T(1)}.\]
Its image $P(L)$ is the space $\ell_A$, isomorphic to $A$ as an algebra, and abelian as a Lie subalgebra of $\mbox{End}(A)$. Combined with any odd, square-zero element $\Delta$ of $\mbox{End}(A)$ (i.e. an inner derivation), we obtain on $A$ the $L_{\infty}$ brackets 
\[ [a_1,\dots,a_n]_{\Delta} =[\cdots [\Delta,\ell_{a_1}],\cdots,\ell_{a_n}]1=\Gamma_{\Delta}^n(a_1,\dots,a_n)1.\;\;\;\Box\]


\subsection{An alternative: Replace subalgebra with cohomology}


\subsubsection{Derived brackets in VOA's and TVOA's}

A {\bf vertex operator algebra (VOA)} $V$ (see e.g.~\cite{Ak2}) comes with an embedding
\begin{equation} \label{voamap} f:V=\oplus_{n\geq n_0}V_{[n]}\rightarrow \mbox{End}(V)[\,[z,z^{-1}]\,],\;\;
 v\mapsto v(z)=\sum_nv_nz^{-n-1}.\end{equation}
A graded VOA is a direct sum $V=\bigoplus_{n\geq n_0}\bigoplus_{m\in\mbox{\bf Z}}     V_{[n]}^{|m|}$ where $[\;]$ denotes the ``weight'' and $|\;|$ denotes the ``super degree''. 
The map $f$ above endows $V$ with products $v\times_nw=v_n(w)$. The formal series $v(z)$ is a {\it vertex operator} and the maps $v_n$ are called the {\it modes} of the vertex operator. 
Besides the {\it standard representation} of $v(z)$ in~(\ref{voamap}), we also have the {\it weight representation} 
$v(z)=\sum_nv_{[n]}z^{-n-[v]}$,
where $v_{[n]}=v_{n-1+[v]}$. The {\it Virasoro operator} ${\cal L}(z)=\sum_{n}{\cal L}_{[n]}z^{-n-2}$ is typically written in this fashion. The special product $\times_{-1}$ is the {\it Wick product} on $V$. 
We call the  graded commutator of $\times_{-1}$ the {\it Wick bracket}.
The {\it residue product} $\times_0$ is a derivation of all products $\times_n$ and 
is left Leibniz on $V$. 
A {\bf topological vertex operator algebra (TVOA)} is a VOA with an odd, square-zero derivation $Q$ (the {\it BRST operator}) of the Wick product and an odd operator $b(z)=\sum_nb_{[n]}z^{-n-2}$ satisfying $[Q,b(z)]={\cal L}(z)$ ($b_{[n]}^2=0$).
\begin{example} With the help of a suitable derivation $d$ such as an odd square-zero residue $d=u_0$, and the residue product $\times_0$, we define a derived bracket 
$ C_d^2(v,w)=(dv)_0w$ 
that is Leibniz on a TVOA $V$. This is how a generalized Batalin-Vilkovisky bracket was defined by Lian and Zuckerman in~\cite{LZ}. Their choice of $d$ was the residue $b_{[-1]}$ of the operator $b(z)$, and the derived bracket was of the form $(-1)^{|v|}(b_{[-1]}v)_0w$. The bracket induces a G-bracket in the BRST cohomology, which is a graded commutative associative algebra induced by the Wick product.\label{lzone}
\end{example}
The BV bracket in \cite{LZ} has attracted attention mostly due to a second and equivalent formulation in the same article as the deviation of $b_{[0]}$ from being a derivation of the Wick product. We first state the following result.
\begin{prop} \cite{Ak3} The modes $v_n$ of a vertex operator $v(z)$ in the standard representation are higher order differential operators with respect to the Wick product as in Definition~\ref{ak}. All $v_n$ with $n\leq -1$ are of order at most zero, and each $v_n$ with $n\geq 0$ is of order at most $n+1$.\label{modes}
\end{prop}
\begin{example} \label{lztwo} Another expression for the BV bracket in Example~\ref{lzone} was obtained in \cite{LZ} by using the order-two differential operator $b_{[0]}$ and the Wick product: the bracket was given by $\pm\Phi_{b_{[0]}}^2$ with respect to Wick. This emulates the BV-type construction, but the bracket is defined on the cohomology.
\end{example}

\subsubsection{Cohomological construction of sh-Lie algebras}

Let us now go beyond the two existing derived bracket constructions which involve restrictions/projections of brackets to an abelian subalgebra and state the main result of this section.
\begin{thm}[Cohomological Derived Brackets Theorem] \label{CDBT}
Let $(V,m,Q)$ be a differential algebra, where $Q$ is an odd, square-zero derivation of $m$, and $\Delta$ be an odd, square-zero operator on $V$ such that $[Q,\Delta]=Q\Delta +\Delta Q=0$. Assume that $m$ is commutative and associative in the cohomology. Then the brackets
\[ \Phi^r_{\Delta}(a_1,\dots,a_r)\]
descend to $H(V,Q)=\mbox{\em Ker}(Q)/\mbox{\em Im}(Q)$ and define an sh-Lie structure on $H(V,Q)$. 
\end{thm}


\begin{remark}
In many examples $m$ is unital pre-Lie, thus there exists a Lie bracket on $V$, which becomes abelian on the cohomology (not a subalgebra).
\end{remark}


\begin{cor}
There exists an $L_{\infty}$ structure on the cohomology $H(V,Q)$ of any differential VOA $(V,\mbox{\em Wick},\mbox{\bf 1},Q)$ with an odd, square-zero operator $\Delta$ of weight zero on $V$ such that $[Q,\Delta]=0$ on $V_{[0]}$ and $\Delta 1=0$. If $\Delta$ is also a mode $v_k$ of some vertex operator ($k\geq 1$), then the $n$-ary brackets after $n=k+1$ vanish.
\end{cor}

{\it Proof of the Cohomological Derived Brackets Theorem.} We use induction on $r$ to show that $\Phi_{\Delta}^r$ preserves both $\mbox{Ker}(Q)$ (a subalgebra of $V$) and $\mbox{Im}(Q)$ (an ideal of $\mbox{Ker}(Q)$). Since $\Phi_{\Delta}^1=\Delta$ anticommutes with $Q$, these subspaces are clearly closed under $\Delta$. Now assume the closure statements hold for $\Phi_{\Delta}^r$. For $a_1,\dots,a_{r+1}\in V$, the $(r+1)$st operator is defined by
\begin{eqnarray} \Phi_{\Delta}^{r+1}(a_1,\dots,a_r,a_{r+1})&=&\Phi_{\Delta}^{r}(a_1,\dots,a_ra_{r+1})\nonumber\\
&&-\Phi_{\Delta}^{r}(a_1,\dots,a_r)a_{r+1}\pm a_r\Phi_{\Delta}^{r}(a_1,\dots,a_{r-1},a_{r+1}).\nonumber\end{eqnarray}
If all $a_i$ are in the kernel, then so is $a_ra_{r+1}$, and the first term on the right-hand side is in $\mbox{Ker}(Q)$ by the induction hypothesis. The remaining two terms consist of a value of the operator $\Phi^r_{\Delta}$ (again in $\mbox{Ker}(Q)$) times an element of the kernel. Finally, assume that all $a_i$ are in $\mbox{Im}(Q)$, and use similar reasoning to complete induction. Since the brackets $\Phi^r_{\Delta}$ are well-defined on the cohomology, we obtain an $L_{\infty}$ algebra on $H(V,Q)$ solely because this space is a  graded commutative associative algebra (Proposition~\ref{voronovex}).~$\Box$

In Part~II of this article, we will strive to give the big picture in sh-Lie constructions and explore relationships with homological perturbation theory (HPT) and deformations (see~\cite{HS}).

{\bf Acknowledgments.} We are grateful to J.~Stasheff for timely and fruitful discussions of the earlier versions of our work. We would also like to thank the anonymous referees and the editor for many historical and technical corrections and for helping us improve the exposition.


F\"{u}sun Akman (akmanf@ilstu.edu), Lucian M. Ionescu (lmiones@ilstu.edu)

\end{document}